\documentclass{article}
\usepackage{graphicx} 

\usepackage[utf8]{inputenc}
\usepackage{graphicx}
\usepackage{amssymb}
\usepackage{amsmath}
\usepackage{graphicx}
\usepackage{enumitem}
\usepackage{amsthm}
\usepackage{tikz-cd}
\newcommand{\R}{\mathbb{R}}
\newcommand{\C}{\mathbb{C}}

\newcommand{\Z}{\mathbb{Z}}
\newcommand{\N}{\mathbb{N}}

\theoremstyle{plain}
\newtheorem{theorem}{Theorem}[section]
\newtheorem{prop}[theorem]{Proposition}

\newtheorem{lemma}[theorem]{Lemma}
\newtheorem{corollary}[theorem]{Corollary}
\theoremstyle{definition}

\newtheorem{remark}[theorem]{Remark}

\renewcommand{\a}{\tilde{a}}
\renewcommand{\b}{\tilde{b}}
\renewcommand{\c}{\tilde{c}}
\renewcommand{\d}{\tilde{d}}
\renewcommand{\S}{\mathcal{S}}

\usepackage{graphicx,trimclip,scalerel}
\def\dblHtextstyle{\clipbox{0pt 0pt 4pt 0pt}{$\mathcal{H}$}\smash{
  \kern-1pt\scalebox{.75}[1.06]{\raisebox{.450pt}{\textsl{\sffamily +}}}\kern-5.15pt
  \kern-1pt\scalebox{.75}[1.04]{\raisebox{.463pt}{\textsl{\sffamily +}}}\kern-5.15pt
  \kern-1pt\scalebox{.75}[1.02]{\raisebox{.476pt}{\textsl{\sffamily +}}}\kern-5.15pt
  \kern-1pt\scalebox{.75}[1.00]{\raisebox{.489pt}{\textsl{\sffamily +}}}\kern-1pt
  \clipbox{4pt -2pt 0pt 0pt}{$\mathcal{H}$}}}
\def\dblH{\scalerel*{\dblHtextstyle}{\mathcal{H}}}
\newcommand{\M}{\mathcal{M}}

\usepackage[
backend=biber,
style=alphabetic,
]{biblatex}
\newcommand*\pFqskip{8mu}
\catcode`,\active
\newcommand*\pFq{\begingroup
        \catcode`\,\active
        \def ,{\mskip\pFqskip\relax}%
        \dopFq
}
\catcode`\,12
\def\dopFq#1#2#3#4#5{%
        {}_{#1}\phi_{#2}\biggl[\genfrac..{0pt}{}{#3}{#4};#5\biggr]%
        \endgroup
}

\title{The Askey-Wilson Function is the Laplace Transform of a Positive Measure}
\author{Aaron Agulnick}
\date{November 2, 2026}

\begin{document}

\maketitle

\begin{abstract}
    The Askey-Wilson function $\phi(x, \gamma)$ is a meromorphic continuation of the Askey-Wilson polynomials in their spectral variable off the lattice $\gamma = \a q^m, m \in \Z_+$. We show that the analytic part of the Askey-Wilson function $\phi^{an}(x, \a q^\gamma)$ is the Laplace transform of a positive measure of total measure $\phi^{an}(x, s_0) < \infty$, under certain conditions on the parameters governing the function. We deduce log-convexity of this function, as well as positive-definiteness and an integral formula for a certain intertwining operator associated with the Askey-Wilson operator $L$. This offers a $q$-analogue of several results, in the trigonometric limit due to Sahi-McSwiggen, and in the rational (i.e. Dunkl) limit due to R\"osler.
\end{abstract}


\section{Introduction}

Integral formulas appear throughout the classical study of special functions to answer questions about asymptotics and other analytic properties, which might be otherwise opaque from the function's origins. For instance, the Dunkl kernel is a special function defined by a joint eigenvalue equations for Dunkl operators, certain differential-reflection operators arising from root systems. There, one demands \begin{equation*}
    T_\xi E_k(x, y) = \langle y, \xi\rangle E_k(x, y)
\end{equation*}
for any direction $\xi$, and obtains the Dunkl kernel as the unique solution $E_k(x,y)$. A parameter specialization sends the Dunkl operators to usual differential operations $\partial_\xi$, and hence the Dunkl kernel generalizes the exponential $e^{xy}$. One can prove in type $A_1$ (e.g. \cite[Example 2.1]{rosler_asymptotic_2002}) the integral formula for the Dunkl kernel \begin{equation*}
    E_k(x,y) = \frac{\Gamma(k+1/2)}{\Gamma(1/2)\Gamma(k)} \int_{-1}^1 e^{txy}(1-t)^{k-1}(1+t)^kdt
\end{equation*}
which, after symmetrization by the Weyl group $S_2 = \{\pm 1\}$, reduces to a classical integral formula for the Bessel function. Hence root systems provide a natural context for vast higher-dimensional generalizations of classical special function properties, including these integral formulas. A key tool in the integral formula above is the \emph{Dunkl intertwiner} $V_k$, which is defined by \begin{equation*}
    T_\xi \cdot V_k(f) = V_k(\partial_\xi f).
\end{equation*}

Recent work of Sahi-McSwiggen \cite{mcswiggen_opdam-cherednik_2026} offers an integral formula for the Opdam-Cherednik kernel, which is the object given analogously to the Dunkl kernel by a joint eigenvalue problem associated to the trigonometric Cherednik algebra. They derive from this integral formula several interesting analytic properties, including the positivity of the Opdam-Cherednik intertwiner, by methods quite dissimilar from the Dunkl setting, where positivity of $V_k$ was proven by R\"osler \cite{rosler_positivity_1999}. \\

We turn our attention to the $q$-deformed setting of the double affine Hecke algebra, for which the trigonometric Cherednik algebra considered in \cite{mcswiggen_opdam-cherednik_2026} is a $q \to 1$ limit. We focus on type $\check{C_1}C_1$, due to independent interest in the algebra and the special functions built from it.  The spherical subalgebra of the DAHA in type $\check{C_1}C_1$ is a $q$-deformation of the coordinate ring of the monodromy manifold of Painlev\'e VI \cite{mazzocco_confluences_2016}. The eigenfunctions of its distinguished generator (acting on the polynomial representation by the Askey-Wilson $q$-difference operator) are a family of orthogonal polynomials known as the Askey-Wilson polynomials, which can be meromorphically continued in the spectral variable to obtain the Askey-Wilson function. The DAHA in type $\check{C}_nC_n$ was first constructed by Sahi in \cite{sahi_nonsymmetric_1999} and its connections to the Askey-Wilson polynomials in rank 1 was explored in \cite{ismail_properties_2000}. The analytic properties of the Askey-Wilson function have been studied extensively (\cite{koelink_askey-wilson_2001}, \cite{stokman_expansion_2002}) due in part to their deep link with the representation theory of $SU(1,1)$ \cite{koelink_fourier_2001}. We hope to continue the study of the Askey-Wilson transform by contributing some analytic properties of its kernel. \\

Our main theorem is that (a certain gauge transformation of) the Askey-Wilson function, an analogue of the Opdam-Cherednik kernel and Dunkl kernel in the trigonometric and rational cases, is also the Laplace transform of a positive measure. More specifically, we prove \begin{theorem}[Theorem \ref{main theorem}]
    Let $\varphi(x, \gamma) := \phi^{an}(x, \a q^\gamma)$ as in (\ref{varphi def}). Suppose the parameters $\alpha$ satisfy the conditions (\ref{Ram-Yip conditions}) and (\ref{coeff conditions}). Then, for each real $x > 0$, $\varphi(x, \gamma)$ is the Laplace transform of a positive measure $\pi_x$ of finite total measure on $\R$.
\end{theorem}

This has numerous pleasant consequences, such as log-convexity, a majorization inequality in the style of \cite{mcswiggen_majorization_2026}, and positivity for a certain DAHA intertwiner. Our techniques differ from \cite{mcswiggen_opdam-cherednik_2026}, as the $q \to 1$ limit is a crucial part of the argument there in moving from a statement about polynomials to a statement about their analytic continuation. We rely on the so-called expansion formula (\ref{expansion formula}) of Stokman \cite{stokman_expansion_2002}, studied in more generality in \cite{stokman_c-function_2014}. Such expansion formulas appear in many similar contexts. For instance, the Jack hypergeometric function introduced in \cite{kaneko_selberg_1993} is defined by such a formula, which is the symmetrization of a formula due to Baker-Forrester \cite{baker_nonsymmetric_1998} expanding the Dunkl kernel in terms of nonsymmetric Jack polynomials. We expect that similar techniques would yield results in those cases, though many of the analytic consequences are already known for the Dunkl kernel; see \cite{rosler_asymptotic_2002}, \cite{rosler_generalized_1998}, and the aforementioned paper \cite{rosler_positivity_1999} on the positivity of the Dunkl intertwiner. The benefit of our approach using the expansion formula is a completely explicit description of the measure (\ref{def pi}) of which $\phi^{an}(x, \a q^\gamma)$ is a Laplace transform; this offers an avenue for more explicit global bounds and a study of the corresponding heat kernel as outlined in \cite[Section 4.2]{mcswiggen_opdam-cherednik_2026}. \\

\subsection{AI/LLM Statement}

All writing and ideas herein are the author's own. Claude Opus 5.5 was used in searching for certain references, in the derivation of the inequalities (\ref{Ram-Yip conditions}) and (\ref{coeff conditions}), and in shortening the proof of Lemma \ref{expansion positivity}; subsequently all work was checked and written solely by the author.

\subsection{Acknowledgments}

The author's research was partially supported by the NSF grants DMS-2601486 and DMS-1645877. The author would like to thank Siddhartha Sahi, Milen Yakimov, and W. Riley Casper for several valuable comments and conversations.

\section{Background}

\subsection{DAHA of Type $\check{C_1} C_1$}

Fix parameters $\alpha = (a,b,c,d) \in \C^4$ and $q \in (0,1)$. We mainly follow the notation of \cite{koornwinder_automorphisms_2025}. The double affine Hecke algebra (DAHA) of type $\check{C_1}C_1$, denoted $\dblH_{a,b,c,d} = \dblH_\alpha$ is the (unital, associative) $\C$-algebra with invertible generators $T_1, T_0, Z$ subject to the relations \begin{equation}
    (T_1+ab)(T_1+1) = 0
\end{equation} 
\begin{equation}
    (T_1Z+a)(T_1Z+b)=0
\end{equation}
\begin{equation}
    (qT_0Z^{-1}+c)(qT_0Z^{-1}+d) = 0
\end{equation}
\begin{equation}
    (T_0+q^{-1}cd)(T_0+1)=0.
\end{equation}

By an \emph{automorphism of DAHAs} we mean an identification $(a,b,c,d) \mapsto (a',b',c',d')$ of parameters and an isomorphism of $\C$-algebras $\dblH_{a,b,c,d} \cong \dblH_{a',b',c',d'}$. We observe from the relations above the following important automorphisms of DAHAs of type $\check{C_1}C_1$ (\cite{koornwinder_automorphisms_2025}): \begin{equation} \label{sigma}
    \sigma: (a,b,c,d) \mapsto (\a, \b, \c, \d)
\end{equation}
\begin{equation*}
    (T_1, T_0, Z) \mapsto (T_1, \sqrt{q^{-1}abcd} T_1^{-1}Z^{-1}, aT_1^{-1}T_0^{-1}).
\end{equation*}
where $\a := \sqrt{q^{-1}abcd}$, with a branch of the square root fixed so that $\sqrt{z} > 0$ when $z > 0$. The other parameters are determined by \begin{equation}
    \label{dual params}
    \a\b = ab, \a\c=ac,\a\d=ad,
\end{equation} 
given explicitly by \begin{equation*}
    \b = \pm \sqrt{qabc^{-1}d^{-1}}, \c = \pm\sqrt{qab^{-1}cd^{-1}}, \d = \pm\sqrt{qab^{-1}c^{-1}d}
\end{equation*}
where the signs are chosen so that (\ref{dual params}) holds. The other key automorphism for us will be
\begin{equation}
    \tau: (a,b,c,d) \mapsto (a,b,c,q/d)
\end{equation}
\begin{equation*}
    (T_1, T_0, Z) \mapsto (T_1, q^{-1}cZT_0^{-1}, Z).
\end{equation*}

For a collection of parameters $\alpha = (a,b,c,d)$, we will write $\alpha_\sigma = \sigma(\alpha) = (\a, \b, \c, \d)$ as in (\ref{sigma}) due to the ubiquity of this automorphism. We call the parameters $\alpha_\sigma$ the \emph{dual} parameters to $\alpha$. Similarly we will write $\alpha_\tau$ for parameters twisted by $\tau$. \\

Let $e$ be the idempotent $e := (1-ab)^{-1}(T_1+1)$. (In the more general presentation of the DAHA, $e$ is the symmetrizer idempotent $\frac{1}{|W|}\sum_{w \in W} w$.) Define the \textbf{spherical subalgebra} $e\dblH_\alpha e$ as the image of $H \mapsto eHe, H \in \dblH_\alpha$. The spherical subalgebra is generated by its elements $De, Xe$ for \begin{equation}
    D = (T_1T_0)+q^{-1}abcd(T_1T_0)^{-1},
\end{equation}
\begin{equation}
    X = Z+Z^{-1}.
\end{equation}

We will need a competing presentation of the DAHA of type $\check{C_1}C_1$ as in \cite{orr_specializations_2018}. The DAHA of (nonreduced) type $\check{C_1}C_1$ is constructed there with parameters $v_{\alpha_1}, v_{2\alpha_1}, v_{\alpha_0}, v_{2\alpha_0}$, the square roots of the parameters in Sahi's presentation \cite{sahi_nonsymmetric_1999} as noted in \cite[Section 2.10, p. 8]{orr_specializations_2018}. From \cite[p. 399]{ismail_properties_2000} or by comparing the presentations we have \begin{equation} \label{parameter correspondence}
    a = v_{\alpha_1}v_{2\alpha_1}, b = -\frac{v_{\alpha_1}}{v_{2\alpha_1}}, c = q^{1/2}v_{\alpha_0}v_{2\alpha_0}, d = -q^{1/2}\frac{v_{\alpha_0}}{v_{2\alpha_0}}.
\end{equation}



\subsection{Askey-Wilson Polynomials and Function}


We make the following standing assumptions on the DAHA parameters: $q \neq 0$, $q$ not a root of unity, $a, b, c, d \neq 0$, $abcd \neq q^{-m}$ for any $m \in \N$, and $\mathrm{Re}(a) > 0$. \\

The DAHA $\dblH_\alpha$ acts on the space $\M$ of Laurent polynomials $f \in \C[z, z^{-1}], f(z) = f(z^{-1})$ by \begin{equation*}
    T_0 \cdot f = \frac{(q^{-1}cd+1)z^2-(c+d)z}{q-z^2}f(z) - \frac{(c-z)(d-z)}{q-z^2}f(qz^{-1}),
\end{equation*}
\begin{equation*}
    (T_1 \cdot f)(z) = \frac{(a+b)z-(1+ab)}{1-z^2}f(z) + \frac{(1-az)(1-bz)}{1-z^2}f(z^{-1}),
\end{equation*}
\begin{equation*}
    (Z \cdot f)(z) = zf(z).
\end{equation*}

Consider also the restriction of this action to $\M^+$, the space of \emph{symmetric} Laurent polynomials. Both representations are faithful \cite[Theorem 3.2]{sahi_nonsymmetric_1999}. Observe also that $De \in e\dblH_\alpha e$ acts on $f \in \M^+$ by $L = L_{a,b,c,d}$ the \textbf{Askey-Wilson operator} $$Lf(z) := (1+q^{-1}abcd)f(z)$$ 
$$+ \frac{(1-az)(1-bz)(1-cz)(1-dz)}{(1-z^2)(1-qz^2)}(f(qz)-f(z))$$ 
$$+ \frac{(a-z)(b-z)(c-z)(d-z)}{(1-z^2)(q-z^2)}(f(q^{-1}z)-f(z))$$
and that $(Xe \cdot f)(z) = (z+z^{-1})f(z)$, see \cite[Section 1]{ismail_properties_2000} and \cite{sahi_nonsymmetric_1999}, and \cite[Proposition 4.11]{koornwinder_automorphisms_2025} for the notation we use here.

\begin{remark}
    In fact, the operators $L$ and $(z+z^{-1})$ viewed in $\mathrm{End}(\M^+)$ generate an algebra $\mathrm{AW}(3, Q_0)$, the quotient of the Zhedanov algebra $\mathrm{AW}(3)$ by the centrally generated ideal $(\Omega - Q_0)$, for $\Omega$ the Casimir and $Q_0$ a constant \cite{zhedanov_hidden_1991}.
\end{remark}

We define the \textbf{symmetric Askey-Wilson polynomials} $R_m(z)$ as the eigenfunctions in $\M^+$ of the Askey-Wilson operator $L$. The eigenvalue of $R_m(z)$ is \begin{equation} \label{eigenval eqn}
    L \cdot R_m(z) = \lambda_m R_m(z), \hspace{1cm} \lambda_m := q^{-m} + abcdq^{m-1}.
\end{equation}

Under our standing assumptions on the DAHA parameters, all the eigenvalues $\lambda_m$ are distinct. \\

\begin{prop}
    We have
    \begin{equation} \label{hypergeom def}
        R_n(z) = \pFq{4}{3}{q^{-n}, q^{n-1}abcd,az,az^{-1}}{ab, ac, ad}{q, q}
    \end{equation}
\end{prop}
where $_r\phi_s$ is the basic hypergeometric series as in e.g. \cite[(1.10.1)]{koekoek_basic_2010}. \\

We define the \textbf{Askey-Wilson function} $\phi(x, \gamma)$ as a kind of analytic continuation of $R_n(z)$ in the spectral parameter: namely, we set \begin{equation}
    \phi(x, \gamma) := \frac{(qax\gamma/\d, qa\gamma/\d x, qabc/d;q)_\infty}{(\a\b\c \gamma, q\gamma/\d, qx/d, q/dx, bc, qb/d, qc/d; q)_\infty} \: _8W_7(\a \b \c \gamma/q; ax, a/x, \a \gamma, \b \gamma, \c \gamma; q, q/\d \gamma)
\end{equation}
for $|q/\d \gamma| < 1$, where $(x; q)_\infty$ is the $q$-Pochhammer symbol \begin{equation*}
    (x;q)_\infty = \prod_{k = 0}^\infty (1-xq^k), \hspace{1cm} (x,q)_n = \prod_{k=0}^{n-1}(1-xq^k).
\end{equation*}

Let $s_m = \a q^{n}$ and denote by $\S = \{s_m: m \in \Z+\}$, so that \begin{equation} \label{analytic continuation}
    \phi(z, s_m) = R_m(z)
\end{equation}
for $m \in \Z_+$. We also define the \textbf{Gaussian} \begin{equation}
    G(x) := (dx, d/x;q)_\infty^{-1}.
\end{equation}

\begin{remark}
    The name ``Gaussian" goes back to Cherednik (\cite{cherednik_one-dimensional_2001}, \cite{cherednik_difference_1997}) in connection with Macdonald polynomials; it controls the analytic properties of the hypergeometric series $\phi$, see below. A key property is that conjugation by the Gaussian implements the automorphism $\tau$ on the standard representation of $\dblH_\alpha$, see \cite[Eqs. (4.7)-(4.8)]{koornwinder_automorphisms_2025}. This has analogues in several other contexts, such as the $q$-Toda system and quantum Q-system \cite{di_francesco_macdonald_2024}. \\
\end{remark}

Finally, define the \textbf{analytic part of the Askey-Wilson function} \begin{equation} \label{def phi^an}
    \phi^{an}(x, \gamma) = G^\tau(x)^{-1}G^{\sigma \tau}(\gamma)^{-1}\phi(x, \gamma).
\end{equation}

For convenience, we set \begin{equation}
    \label{varphi def}
    \varphi(x, \gamma) := \phi^{an}(x, \a q^\gamma)
\end{equation} so that $\varphi(z, m) = G^\tau(z)^{-1}G^{\sigma\tau}(s_m)^{-1}R_m(z)$. Not only is this a natural way to view the Askey-Wilson function as an analytic continuation of the Askey-Wilson polynomials off of the \emph{integer} lattice, but it is also technically convenient in analogy with the Macdonald polynomial convention of \cite{mcswiggen_opdam-cherednik_2026}, where a monomial $x^\nu$ is evaluated by $x \mapsto q^{\langle x, \nu \rangle}$. The point is to view the dominant weights as $\S = \{\a q^m: m \in \Z_+\}$ instead of $\Z_+$ itself. \\

Following \cite{stokman_expansion_2002}, we write $R^\sigma_m$ or $\widetilde{R}_m$ for the Askey-Wilson polynomial with \emph{dual} parameters $\alpha_\sigma = (\a, \b, \c, \d)$; equivalently, $R^\sigma_m$ is the solution of the eigenvalue problem (\ref{eigenval eqn}) in the polynomial representation of $\dblH_{\alpha_\sigma}$. In general, if $H$ depends on the DAHA parameters $\alpha$, and if $\eta$ is an automorphism of DAHAs of type $\check{C_1}C_1$, we write $H^\eta$ for the object $H$ with parameters $\eta(\alpha)$. \\

\begin{theorem}[\cite{stokman_expansion_2002}, Theorem 4.2] 
    The analytic part of the Askey-Wilson function $\phi^{an}(x, \gamma)$ admits a biorthogonal expansion in terms of the Askey-Wilson polynomials, namely \begin{equation} \label{expansion formula}
        \phi^{an}(x, \gamma) = G^{\tau \sigma \tau}(s_0^\tau)^{-1} \sum_{m \in \Z_+} h^{\tau \sigma}(s_m)G^{\tau \sigma \tau}(s_m)R_m^{\sigma \tau}(\gamma)R_m^\tau(x).
    \end{equation}
    with \begin{equation*}
        h^{\tau \sigma}(s_m^\tau) = \frac{(1-abcq^{2m}/d)(abc/d, ab, ac, qa/d;q)_m}{(1-abc/d)(q, bc, qb/d, qc/d;q)_m}a^{-2m}.
    \end{equation*}
\end{theorem}

From the theorem above, it is obvious that \begin{equation}
    \phi^{an}(x, \gamma) = \widetilde{\phi^{an}}(\gamma, x)
\end{equation}
and hence \begin{equation}
    \phi(x, \gamma) = \widetilde{\phi}(x, \gamma).
\end{equation}

As a consequence of (\ref{analytic continuation}), we see also that \begin{equation}
    R_m(s_n) = \widetilde{R_n}(s_m)
\end{equation}
or equivalently by $z \mapsto z^{-1}$ symmetry \begin{equation*}
    R_m(s_n^{-1}) = \widetilde{R_n}(s_m^{-1}),
\end{equation*}
the form preferred in \cite{koornwinder_dualities_2018}. \\

These three facts are known as duality for the Askey-Wilson polynomials and function, originating from Sahi's work \cite{sahi_nonsymmetric_1999} and \cite{ismail_properties_2000}. See also the treatment in \cite{koornwinder_dualities_2018}.

\section{Main Results}

This section is essentially dedicated to the proof of our main theorem, Theorem \ref{main theorem}, which will rely on two key positivity results which we offer first as lemmas.

\begin{lemma} \label{Ram-Yip coeff positivity}
    Suppose \begin{equation}
        \label{Ram-Yip conditions}
        a > 0 > b, c > 0 > d, a < |b|, a|b|<1, c< |d|, c|d|<q.
    \end{equation}
    Then, writing the Askey-Wilson polynomial as \begin{equation} \label{poly decomp}
        R_m(x) = c_m(x^m+x^{-m}) + \sum_{n < m} c_n(x^n+x^{-n}),
    \end{equation}
    all the coefficients $c_n \geq 0$ and $c_m > 0$.
\end{lemma}
\begin{proof}
The eigenfunction equation (\ref{eigenval eqn}) for $R_m$ and \cite[Proposition 2.21]{orr_specializations_2018} imply that $R_m$ are proportional to Orr-Shimozono's symmetric Macdonald-Koornwinder polynomials in type $\check{C_1}C_1$. Furthermore, Orr-Shimozono's convention is to normalize so the top-degree (in $x^n+x^{-n}$) coefficient is 1, while our convention (following Koornwinder-Mazzocco) has top-degree coefficient \begin{equation*}
    \frac{a^n(q^{n-1}abcd;q)_n}{(ab,ac,ad;q)_n}
\end{equation*}
from analyzing (\ref{hypergeom def}) which is positive under the assumptions (\ref{Ram-Yip conditions}). So it suffices to use the Ram-Yip formula \cite[Proposition 3.20]{orr_specializations_2018} to check positivity of the coefficients $c_n$. \\ 

The coefficients are positive when \begin{equation}
        \label{OS parameter positivity}
        (v_{\alpha_1}, v_{2\alpha_1}, v_{\alpha_0}, v_{2\alpha_0}, q) \in (0,1)^5.
    \end{equation}
    Indeed, the scalars $c_j, d_j$ appearing in $f_J(u)$ of that proposition are differences $v_\alpha^{-1}-v_\beta$ for $\alpha,\beta \in \{\alpha_1, 2\alpha_1, \alpha_0, 2\alpha_0\}$, so certainly are positive when all the $0 < v_\alpha < 1$. By \cite[Remark 3.17]{orr_specializations_2018}, the $\xi_j$ are products of a positive power of $q$ and positive powers of the Hecke parameters $v_{\alpha_i}, v_{2\alpha_i}$; so when all the $v_{\alpha_i}, v_{2\alpha_i} < 1$, certainly $\xi_j < 1$, and so the denominators $1-\xi_j^2 > 0$. Likewise the $v_{uu_\lambda^{-1}}$ and $v_{\mathrm{dir}(p_J(u))}$ are products of Hecke parameters and their inverses, which are thus positive. \\

    Then, using the parameter correspondence (\ref{parameter correspondence}), the conditions (\ref{OS parameter positivity}) are equivalent to (\ref{Ram-Yip conditions}).
\end{proof}
\begin{remark}
    In \cite{mcswiggen_opdam-cherednik_2026}, the parameters $v_\alpha$ are assigned by $v_\alpha = q^{k_\alpha}$, $k \geq 0, q \in (0, 1)$ for $k_\alpha$ a class function, i.e. a Weyl-group-invariant function. Hence our assumptions (\ref{Ram-Yip conditions}) are equivalent to those made by Sahi-McSwiggen in the particular case of type $\check{C_1}C_1$, except that we do not assume $\alpha_\sigma = \alpha$.
\end{remark}

\begin{lemma} \label{expansion positivity}
    Suppose \begin{equation}
        \label{coeff conditions}
        a, c \in (0, 1); \hspace{1cm} b, d < 0; \hspace{1cm} q|b|<|d|; \hspace{1cm} abc/d < 1.
    \end{equation}
    Then, for fixed real $z > 0$ and every $s_m \in \S$, \begin{equation}
        G^{\tau \sigma \tau}(s_0^\tau)^{-1}h^{\tau \sigma}(s_m)G^{\tau \sigma \tau}(s_m)R_m^\tau(z) \geq 0.
    \end{equation}
\end{lemma}
\begin{proof}
    By the identities \begin{equation}
        h^{\tau \sigma}(s_m^\tau) = \frac{(1-abcq^{2m}/d)(abc/d, ab, ac, qa/d;q)_m}{(1-abc/d)(q, bc, qb/d, qc/d;q)_m}a^{-2m}
    \end{equation}
    and \begin{equation}
        G^{\tau \sigma \tau}(s_0^\tau)^{-1}G^{\tau \sigma \tau}(s_m) = \frac{(bc;q)_m}{(qa/d;q)_m}(\frac{-a}{d})^m q^{m(m+1)/2}
    \end{equation}
    (see \cite[p. 317]{stokman_expansion_2002}), we see that $h^{\tau \sigma}(s_m^\tau)$ and $G^{\tau \sigma \tau}(s_0^\tau)^{-1}G^{\tau \sigma \tau}(s_m)$ are nonnegative for parameters satisfying (\ref{coeff conditions}). Hence it remains to check that $R_m^\tau(x) \geq 0$ there. \\

    Use $x := \frac{1}{2}(z+z^{-1})$ and observe that $x > 0$ iff $z > 0$. The Askey-Wilson polynomial $R_m(z)$ depends only on $z+z^{-1}$, i.e. only on $x$, so write $p_n(x; \alpha) = a^{-n}(ab, ac, ad;q)_nR_n(z)$. This is the usual normalization of the Askey-Wilson polynomials, in particular symmetric in permutations of parameters $(a,b,c,d)$. The generating function \cite[(14.1.14)]{koekoek_basic_2010}, applied with parameters $\alpha_\tau = (a,b,c,q/d)$ and $z = e^{i\theta}, x = \cos \theta = \frac{1}{2}(z+z^{-1})$ has $t^n$-coefficient \begin{equation}\label{gen func}
        \frac{p_n(x; a,b,c,q/d)}{(ac, bq/d, q;q)_n} = \sum_{k=0}^n \frac{(bz, (q/d)z;q)_k}{(bq/d, q;q)_k}z^{-k} \cdot \frac{(a/z, c/z; q)_{n-k}}{(ac, q; q)_{n-k}}z^{n-k}. 
    \end{equation}
    The numerators of the $z^{-k}$-term are all positive since $b, d,$ and hence $q/d$ are all negative, and since $z > 0$ and $bq/d < 1$. The numerators of the $z^{n-k}$-coefficients are positive when $z > a$ and $z > c$; we will remove this assumption shortly. The denominators $(bq/d; q)_k > 0$ when $bq/d < 1$, and $(ac;q)_{n-k} > 0$ when $ac < 1$, which both follow from our parameter assumptions. Lastly, the denominator of the left-hand side is also positive under these conditions. Hence $p_n(x; \alpha_\tau) = p_n(\frac{1}{2}(z+z^{-1}); \alpha_\tau) > 0$ whenever $z > \max(a,c)$. But $p_n$ is invariant under $z \mapsto z^{-1}$ (equivalently, use that $p_n$ is invariant under parameter permutations, and swap $a \leftrightarrow c, b \leftrightarrow q/d$ in the generating function identity), so applying the same argument to $p_n(\frac{1}{2}(z^{-1}+z); \alpha_\tau)$ shows that $p_n(x; \alpha_\tau) > 0$ when $z < \frac{1}{\max(a,c)}$. When $a,c \in (0, 1)$, this proves $p_n(x; \alpha_\tau)$ for all $z > 0$, and hence all $x > 0$. But then the normalization factor $(a^{-n}(ab, ac, ad;q)_n)^{-1}$ is also positive, so $R_n^\tau(z) > 0$ for all $z > 0$ too.
\end{proof}

Now we are ready to prove the main theorem.

\begin{theorem} \label{main theorem}
    Let $\varphi(x, \gamma) := \phi^{an}(x, \a q^\gamma)$ as in (\ref{varphi def}). Suppose $\alpha$ satisfies the conditions (\ref{Ram-Yip conditions}) and (\ref{coeff conditions}). Then, for each real $x > 0$, $\varphi(x, \gamma)$ is the Laplace transform of a positive measure $\pi_x$ of finite total measure on $\R$.
\end{theorem}
\begin{proof}
    Write $C_m(x) := G^{\tau \sigma \tau}(s_0^\tau)^{-1}h^{\tau \sigma}(s_m)G^{\tau \sigma \tau}(s_m)R_m^\tau(x)$, so that by the expansion formula (\ref{expansion formula}), we have \begin{equation*}
        \phi^{an}(x, \gamma) = \sum_{m \geq 1} C_m(x) R^{\sigma \tau}_m(\gamma).
    \end{equation*}
    Expand $R^{\sigma\tau}_s(\gamma)$ in monomials according to (\ref{poly decomp}) to obtain \begin{equation*} 
        \phi^{an}(x, \gamma) = \sum_{m \geq 1} \sum_{n \leq m} C_m(x) c_n (\gamma^n+\gamma^{-n})
    \end{equation*}

    with the $c_n \geq 0$ by Lemma \ref{Ram-Yip coeff positivity}, as $\alpha_{\sigma\tau}$ satisfies the hypothesis if and only if $\alpha$ does. This is best seen by the fact that $\sigma$ and $\tau$ act on Orr-Shimozono's parameterization by permuting the $v_\alpha$, and the region (\ref{OS parameter positivity}) is invariant under such permutations. \\
    
    In particular, \begin{equation} \label{summation}
        \varphi(x, \gamma) = \phi^{an}(x, \a q^\gamma) = \sum_{m \geq 1} \sum_{n \leq m} C_m(x) c_n (\a^n q^{\gamma \cdot n}+\a^{-n}q^{-\gamma \cdot n}).
    \end{equation}
    Observe that the summation (\ref{summation}) expresses $\varphi(x, \gamma)$ as the Laplace transform of the measure \begin{equation} \label{def pi}
        \pi_x := \sum_{m \geq 1} \sum_{n \leq m}C_m(x)c_n  (\a^n \delta_{n \log q} + \a^{-n}\delta_{-n\log q})
    \end{equation}
    since \begin{equation} \label{positive combination}
    \int_\R e^{\gamma \cdot y} \pi_x(dy) = \sum_{m \geq 1} \sum_{n \leq m} C_m(x)c_n(\a^n q^{\gamma \cdot n}+\a^{-n}q^{-\gamma \cdot n}) = \varphi(x, \gamma).
\end{equation}
    By the given conditions on the DAHA parameters $a,b,c,d,$ and by Lemmas \ref{Ram-Yip coeff positivity} and \ref{expansion positivity}, the measure $\pi_x$ is positive as desired. Furthermore, the total measure is \begin{equation}
        \sum_{m \geq 1} \sum_{n \leq m} C_m(x)c_n(\a^n+\a^{-n}) = \varphi(x, 0) = \phi^{an}(x, s_0) < \infty
    \end{equation}
    since $\phi^{an}(x, \gamma)$ is analytic by design, so that $\pi_x$ is a finite measure for any $x > 0$.
\end{proof}


\section{Applications}

We turn now to several applications of the main theorem, mostly following \cite{mcswiggen_opdam-cherednik_2026}. The point is that the integral formula (\ref{positive combination}) allows one to simply read off many analytic consequences via properties of the exponential. 

\begin{corollary}
    Under the parameter assumptions (\ref{Ram-Yip conditions}) and (\ref{coeff conditions}), $\varphi(x, \gamma)$ is log-convex in the spectral parameter. That is, the function $\gamma \mapsto \varphi(x, \gamma)$, and hence $m \mapsto G^\tau(x)^{-1}G^{\sigma\tau}(s_m)^{-1}R_m(x)$, is log-convex for any fixed $x > 0$.
\end{corollary}
\begin{proof}
    Formula (\ref{summation}) witnesses $\varphi(x, \gamma)$ as a positive combination of log-convex functions, hence log-convex. The statement about the Askey-Wilson polynomials follows by (\ref{analytic continuation}) and (\ref{varphi def}).
\end{proof}

Furthermore, it is immediate from the positivity conditions in Lemma \ref{Ram-Yip coeff positivity} and the argument of \cite[Proposition 4.3]{mcswiggen_opdam-cherednik_2026} that $m \mapsto R_m(s_n)$ is log-convex in the spectral parameter $m$ whenever the \emph{dual parameters} $\alpha_\sigma$ satisfy conditions (\ref{Ram-Yip conditions}); the parameters need not be self-dual. Indeed, since $\sigma$ acts by permuting the parameters in Sahi's parameterization, and the region (\ref{Ram-Yip conditions}) corresponds to $(0,1)^5$ in that parameterization, log-convexity holds whenever $\alpha$ itself satisfies (\ref{Ram-Yip conditions}).

\begin{corollary}
    Under the same parameter assumptions (\ref{Ram-Yip conditions}) and (\ref{coeff conditions}), the function $\phi^{an}(x, \gamma)$ is nondecreasing in $|\log \gamma|$.
\end{corollary}
\begin{proof}
    By the same argument as in the theorem, but removing the $\a^{\pm n}$ from the measure, we see that $\gamma \mapsto \phi^{an}(x, q^\gamma)$ is the Laplace transform of a positive measure as well. But it is also \emph{even} in $\gamma$, because $\phi^{an}(x, \gamma) = \phi^{an}(x, \gamma^{-1})$, as can be read off the expansion formula (\ref{expansion formula}). Hence $\gamma \mapsto \phi^{an}(x, q^\gamma)$ is nondecreasing in $|\gamma|$, so that $\phi^{an}(x, \gamma)$ is nondecreasing in $|\log \gamma|$.
\end{proof}
\begin{remark}
    The final implication in the preceding corollary can be derived ``by hand" in rank 1: by inversion-invariance, take $\phi^{an}(x, q^\gamma) = \frac{1}{2}(\phi^{an}(x, q^\gamma) + \phi^{an}(x, q^{-\gamma}))$, and then the Laplace transform expansion of the right-hand side shows that the expression is a positive combination of $\cosh(y \cdot \gamma)$, so manifestly nondecreasing in $|\gamma|$. But this is really the majorization principle of \cite{mcswiggen_majorization_2026}. Here $W_0 = \{\pm 1\}$ acts by negation on $\R$, so convexity (from our main theorem) and $W_0$-invariance in $q^\gamma$ (from the expansion formula) give Schur-convexity by \cite[Proposition C.2, p. 97]{marshall_schur-convex_2011}. In particular, $\phi^{an}(x, q^{\gamma}) \geq \phi^{an}(x, q^{\gamma'})$ whenever $\gamma$ is in the convex hull of the orbit $W_0 \cdot \gamma'$, i.e. when $|\gamma| \leq |\gamma'|$. (Thus $\phi^{an}(x, q^\gamma)$ is ``nondecreasing in majorization.") \\
\end{remark}

In connection with \cite{rosler_positivity_1999} and \cite{mcswiggen_opdam-cherednik_2026}, we offer one more idea. For $x > 0$, define the \textbf{intertwiner} $V_\alpha(f)(x) := \int_\R f(y)\pi_x(dy)$ and observe that \begin{equation}
    \label{intertwiner on exponential}
    V_\alpha(e^{-\cdot\gamma})(x) = \int_\R e^{y\cdot \gamma} \pi_x(dy) = \varphi(x, \gamma).
\end{equation}

The function space upon which this calculation occurs is slightly subtle: $V_\alpha(f)$ makes sense for polynomials $f$ for growth rate reasons, since it makes sense for exponentials by (\ref{positive combination}). In order to suitably extend the space of polynomials to a function space containing the exponential, we complete the space of polynomials with respect to the norm $||p||_r = \sum_{n=0}^\infty ||p_n||_{\infty, B_r}$, where $p_n$ is the homogeneous degree-$n$ part of $p$ and $B_r$ is the closed $r$-ball centered at the origin. The completed space is denoted $A_r$, and $e^y \in A_r$ for any $r$. For a complete treatment in the Dunkl setting, see \cite[Section 2.4]{rosler_dunkl_2003}. \\

The property (\ref{intertwiner on exponential}) characterizes a unique continuous linear map $V_\alpha: A_r \to A_r$, which is thus necessarily the intertwiner as defined above. Then, since $\phi(x, \gamma)$ satisfies \begin{equation*}
    L^{(x)} \cdot \phi(x, \gamma) = (-1-\a + \a(\gamma + \gamma^{-1}))\phi(x, \gamma)
\end{equation*}
(see e.g. \cite[Section 3]{koelink_askey-wilson_2001}) and in view of the definitions (\ref{def phi^an}) and (\ref{varphi def}), we have
\begin{equation*}
    G^\tau(x)^{-1} L^{(x)} G^\tau(x) \cdot V_\alpha(e^{-\cdot \gamma})(x) = G^\tau(x)^{-1} L^{(x)} G^\tau(x) \cdot \varphi(x, \gamma)
\end{equation*}
\begin{equation*}
     = (-1-\a^2+\a(\a q^\gamma + \a^{-1}q^{-\gamma})) \varphi(x, \gamma) = \int_\R (-1-\a^2+\a(\a q^\gamma + \a^{-1}q^{-\gamma}))e^{y \cdot \gamma}\pi_x(dy)
\end{equation*}
\begin{equation*}
    = V_\alpha(M \cdot e^{-\cdot \gamma})(x)
\end{equation*}
where $M$ is the operator $f(y) \mapsto (-1-\a^2)f(y) + \a^2 f(y + \log q) + f(y - \log q)$. That is, the operator $V_\alpha$ intertwines the action of $G^\tau(x)^{-1}L^{(x)}G^\tau(x)$ with that of $M$, and by density of exponential polynomials this property equivalently characterizes $V_\alpha$. Then from Theorem \ref{main theorem} we have: \begin{corollary}
    The intertwiner $V_\alpha$ is positivity-preserving on the half-line $\R_+$ whenever the parameters $\alpha$ satisfy (\ref{Ram-Yip conditions}) and (\ref{coeff conditions}), i.e. if $f(x) \geq 0$ for all $x$ then $V_\alpha(f)(x) \geq 0$ for all $x > 0$.
\end{corollary}




\printbibliography

\end{document}